\documentclass%[12pt]
{amsart}

\usepackage{amssymb,amscd}
\usepackage[all,line,arc,curve,color,frame,pdf]{xy}
\usepackage{tikz}
\usepackage{tikz-cd}
\usetikzlibrary{positioning, trees, snakes}
\usepackage[textsize=tiny]{todonotes}
\usepackage[colorlinks = true, citecolor = blue]{hyperref}
\usepackage{enumitem}
\usepackage{float}
\usepackage{ytableau}
\usepackage{algpseudocode,mathrsfs}
\usepackage[symbol]{footmisc}
\numberwithin{equation}{section}
\theoremstyle{plain}
\newtheorem{thm}{Theorem}[section]

\newtheorem{question}[thm]{Question}

\theoremstyle{remark}
\newtheorem{definition}[thm]{Definition}
\newtheorem{assumption}[thm]{Assumption}

\newtheorem{ex}[thm]{Example}

\theoremstyle{definition}

\usepackage[textsize=tiny]{todonotes}

\usepackage{enumitem}

\newcommand\ignore[1]{}

\newcommand\CC{{\mathbb{C}}}
\newcommand\NN{{\mathbb{N}}}
\newcommand\PP{{\mathbb{P}}}

\newcommand\ZZ{{\mathbb{Z}}}

\newcommand{\aS}{\mathfrak{S}}

\def\operatorname#1{\mathop{\rm #1}\nolimits}

\def\codim{\operatorname{codim}}

\def\deg{\operatorname{deg}}

\def\det{\operatorname{det}}

\def\GL{\operatorname{GL}}

\newcommand{\pb}{\ar@{}[dr]|{\text{\pigpenfont J}}}

\makeatletter
\newcommand{\xleftrightarrow}[2][]{\ext@arrow 3359\leftrightarrowfill@{#1}{#2}}
\makeatother
\newcommand{\xdasharrow}[2][->]{
\tikz[baseline=-\the\dimexpr\fontdimen22\textfont2\relax]{
\node[anchor=south,font=\scriptsize, inner ysep=1.5pt,outer xsep=2.2pt](x){#2};
\draw[shorten <=3.4pt,shorten >=3.4pt,dashed,#1](x.south west)--(x.south east);
}}

\title{Equations for $\GL$ invariant families of polynomials}

\author[Breiding]{Paul Breiding}
\address{Max Planck Institute for Mathematics in the Sciences,
Leipzig}
\author[Ikenmeyer]{Christian Ikenmeyer}
\address{University of Liverpool, United Kingdom}
\author[Micha{\l}ek]{Mateusz Micha{\l}ek}
\address{Max Planck Institute for Mathematics in the Sciences,
Leipzig, Germany and University of Konstanz
}
\author[Hodges]{Reuven Hodges}
\address{U.~Illinois at Urbana-Champaign, Urbana, IL 61801, USA}

\thanks{
P.~Breiding has received funding from the European Research Council (ERC) under the European Union's Horizon 2020 research and innovation programme (grant agreement No 787840) and funding from the Deutsche Forschungsgemeinschaft (DFG, German Research Foundation) -- Projektnummer 445466444. C.~Ikenmeyer was supported by the DFG grant IK 116/2-1.}

\begin{document}

\begin{abstract}
We provide an algorithm that takes as an input a given parametric family of
homogeneous polynomials, which is invariant under the action of the general linear group, and an integer $d$. It outputs the ideal of that family intersected with the space of homogeneous polynomials of degree $d$.  Our motivation comes from  Question~7 in \cite{ranestad2020twenty} and Problem~13 in \cite{BerndList}, which ask to find equations for varieties of cubic and quartic symmetroids.
The algorithm relies on a database of specific Young tableaux and highest weight polynomials. We provide the database and the implementation of the database construction algorithm. Moreover, we provide a julia implementation to run the algorithm using the database, so that more varieties of homogeneous polynomials can easily be treated in the future.
%In addition, we conduct a numerical experiment seeking to determine the degree of the variety of quartic symmetroids.
\end{abstract}

\maketitle

\section{Introduction}
Many mathematical models are defined by nonlinear maps $f:V\rightarrow W$ between vector spaces. For instance, such models are common in statistics and physics.  The setting allows to generate possible outcomes of the model, by evaluating $f$. This is called the \emph{forward problem}. On the other hand, the \emph{inverse problem} is to decide if a point $w\in W$ belongs to the image of $f$, and if so, to determine its preimage.

In this article we focus on the case when $f$ is a polynomial map. Under this assumption the forward problem consists in evaluating a system of polynomials, and the inverse problem is to solve a system of polynomial equations.

Our main aim is to describe the closure of the image of $f$, when $V$ and $W$ are complex vector spaces. The goal is to describe the polynomial equations that vanish on the image of $f$. Having such equations at hand decouples the inverse problem: for the decision problem, whether or not $w$ is in the image of $f$, one can evaluate the polynomials at $w$ instead of solving a system of equations. The former is much simpler than the latter.

The classical method to find equations relies on the computation of a lexicographic Gr\"obner basis \cite{cox2013ideals, MiSt} to perform elimination of variables. It is a symbolic method, that in practice may be used only on small examples. Thus, the motivation for us is to describe an alternative algorithm that can go beyond these small cases. In general, this task is too ambitious. But, if we assume that the problem has some underlying symmetries, we can use the power of \emph{representation theory} to reduce complexity. In this paper, we make the following assumption for  $f$: we require it to be mapping into a vector space of polynomials and we assume that the image of~$f$ is GL-invariant.

\begin{assumption}\label{assumption1}
We assume that $W=S^c(\CC^n)$ is the space of homogeneous polynomials of degree $c$ in $n$ variables. Furthermore, we assume that the image of $f$ is invariant under $\GL(\CC^n)$, which acts by variable substitution.
\end{assumption}

Our motivation comes from Question~7 in \cite{ranestad2020twenty} and Problem~13 in \cite{BerndList}, which ask to find equations for varieties of cubic and quartic symmetroids. These are subvarieties of the vector space of homogeneous polynomials in $n=4$ variables of degree respectively $c=3$ and $c=4$. They are   $\mathrm{GL}(4)$-invariant. We address these problems in Section \ref{sec:exm}.

We note that the $\GL$ action both gives us many advantages and is very natural. Our ambient space $S^c(\CC^n)$ of polynomials may be regarded as a space of varieties. Following Felix Klein's Erlangen program \cite{Klein1893} geometric quantities should be group-invariant. Thus, very often when studying sets of polynomials, we would like those sets not to depend on the choice of the basis. This is precisely the $\GL$ invariance. Further, the space of polynomials vanishing is often huge, but the $\GL$ action reduces the complexity and allows us to describe it using just a few generators.

\section*{Acknowledgements}
We thank Bernd Sturmfels for motivating questions. These questions not only inspired us, but also revealed bottlenecks of our algorithms. We also thank anonymous referees, whose suggestions helped us to improve the exposition of the paper.
\section{Contributions}

We present an algorithm to study the image of $f$ under Assumption \ref{assumption1}. This algorithm produces the following: let
$$X:=\overline{\mathrm{im}(f)}$$
be the closure of the image of $f$ and let $I$ be the ideal of polynomials that vanish on~$X$. Given $f$ and any $d$ we return the minimal set of polynomials that under the $\GL$ action span $I_d$. This algorithm is exact, ie.~ does not rely on any approximations.
However, instead of a purely symbolic algorithm that works with the parametrized variety $X$ directly, a much more efficient implementation just samples from $X$ (without approximations) and uses only the sampled points as input,
which reduces the finding problem to a linear algebra problem. The details are given in Section~\ref{sec:algorithm}.
This means that due to unlucky sampling in principle the algorithm could output equations that are not actually equations. In practice the probability of this is extremely low and can be further reduced to an inverse exponentially small probability by running the algorithm several times. Further, a posteriori, it is easy to check if the equations actually vanish on $\mathrm{im}(f)$.
One of the algorithm's central ingredients is a \emph{database} which contains bases of highest weight spaces for different plethysms.

The variety of quartic (resp.~cubic) symmetroids consists of polynomials that are determinants of symmetric four by four (resp.~three by three) matrices with entries that are linear forms in four variables. To distinguish these varieties from the general $X$ we will use another symbol:
\begin{align}\label{def_Q}
&Q_3:= \{\det(x_0A_0 + x_1A_1 + x_2A_2 + x_3A_3)\mid A_i\in \CC^{3\times 3}, A_i^T = A_i, 0\leq i\leq 3\},\\
&Q_4:= \{\det(x_0A_0 + x_1A_1 + x_2A_2 + x_3A_3)\mid A_i\in \CC^{4\times 4}, A_i^T = A_i, 0\leq i\leq 3\}.\nonumber
\end{align}
The varieties are $\GL(\CC^4)$ invariant subvarieties of $S^4(\CC^4)$ (resp.\ $S^3(\CC^4)$) of codimension $10$ (resp.\ $4$). We apply our algorithm to this variety, and we obtain the following result.
\begin{thm}\label{thm:main}
There are no equations for $Q_4$ in degrees up to (including) $8$.

There are no equations for $Q_3$ in degrees up to (including) $10$. In degree $11$ the vector space $I(Q_3)_{11}$ is an irreducible representation of dimension $220$ corresponding to the weight $[15,6,6,6]$. The unique highest weight vector has $23824$ many terms. In degree $12$ the vector space $I(Q_3)_{12}$ has four irreducible components corresponding to weights $[15,9,6,6]$, $[16,8,6,6]$, $[17,7,6,6]$ and $[18,6,6,6]$.
\end{thm}

Our second contribution is a numerical algorithm that computes the degree of variets as described above. We applied it to $Q_3$ confirming the known result \cite{Vainsencher2003} that the variety has degree $305$. However, we were not able to compute the degree of $Q_4$, which we leave as a future challenge.

We combine numerical and symbolic methods in our algorithms. Both the numerical and the symbolic algorithm appeared (explicitly or implicitly by using highest weight polynomials as images of symmetrizations over the wreath product) in particular examples before \cite{bci:10, BI:11, BI:13, HIL:13, AIR:14, CIM:15, CHILO:18, BIP:19, DI:19, IK:20, Abd:02, Ott:13, OR:11, Kum:11, BO:11, Rai:13, DHO:14}.
However, to our knowledge, this is the first general implementation and the first one with which it is possible to check for equations of degree 8 on $S^4(\CC^4)$ or find equations of $Q_3$. This is made possible by the use of an idea that we call
\emph{equivariant hash functions}, see Section~\ref{sec:algorithm}.
We provide the source code of our implementation and an easy to use user interface for future researchers to build upon.

We remark that new algorithms for evaluating highest weight polynomials have been developed very recently in \cite{BDI:20}.
No open source implementation of these algorithms is available, but in a special case (see \cite{DI:19}) the running time improvements seem to be of practical importance.

\section{Representation theory}
Representation theory can be  beneficial for large computations. In one line, it allows one to replace a possibly high dimensional irreducible representation, by a one dimensional subspace---the span of the highest weight vector. We briefly recall the relevant concepts for our setting. For more details we refer to \cite{FH, MiSt}.

Every irreducible, polynomial representation $V=V_\lambda$ of $\GL(\CC^n)$ is associated to a Young diagram~$\lambda$ with at most $n$ rows. Fixing the torus $T\subset\GL(\CC^n)$ of diagonal matrices the representation $V$ of $T$ is decomposable $V=\oplus_{\chi\in\ZZ^n} V_\chi$, where $t v=\chi(t) v$ for $v\in V_\chi$ and $\ZZ^n$ is the lattice of characters of the torus $T$. The lexicographically largest $\chi$, say $\chi_0=(l_1,\dots,l_n)$ is called the highest weight of $V$. The Young diagram $\lambda$ has $l_i$ boxes in the $i$-th row. We have $\dim V_{\chi_0}=1$. The unique up to scaling element of $V_{\chi_0}$ is called the \emph{highest weight vector}.
\begin{ex}
Let $V=S^d(\CC^n)$ be the $d$-th symmetric power of $\CC^n$. It is an irreducible representation. The characters of the torus $\chi$ appearing in the representation correspond to $n$-tuples of nonnegative integers summing up to $d$. The highest weight is $(d,0,\dots,0)\in\ZZ^n$. The highest weight vector is $e_1\cdots e_1$. The associated Young diagram is a row with $d$ boxes.
\end{ex}
More generally, for any representation $V$ of $\GL(\CC^n)$ a vector $v\in V$ is called a \emph{highest weight vector} if it is an image of a highest weight vector in some irreducible representation $V_\chi$ under a $\GL(V)$-equivariant map. If $V=\bigoplus_\lambda V_\lambda^{\oplus a_\lambda}$ is the decomposition of $V$, then $a_\lambda$ equals the dimension of the vector space of highest weight vectors in $V$ of weight $\lambda$. Further, any highest weight vector of weight $\lambda$ uniquely determines a subrepresentation $V_\lambda\subset V$. In other words, representation theory allows to replace a possibly large representation $V$ by much smaller spaces of highest weight vectors.

The main observation is that if $X\subset S^c(\CC^n)$ is $\GL(\CC^n)$ invariant, then $I_d$ is a representation of $\GL(\CC^n)$, which is a subrepresentation of $S^d(S^c((\CC^n)^*))$. The representation $S^d(S^c(\CC^n)^*)$ is known as a plethysm. In general the formulas for its decomposition into irreducible representations are not known, and determining a combinatorial description for the multiplicities of irreducibles is Problem 9 in Stanley's list of open problems in algebraic combinatorics \cite{sta:00}.
However, they are known up to $d\leq 5$ \cite{kahle2016plethysm} and for fixed $d$ and $c$ there are algorithms to find such decompositions.
For general $d$ and $c=3$ even the task of deciding positivity of $a_\lambda$ is NP-hard, see \cite{FI:20}.

If $S^d(S^c((\CC^n)^*))=\bigoplus_{\lambda\vdash dc} (S^\lambda)^{\oplus a_\lambda}$ is the decomposition, then we seek to find subrepresentations $(S^\lambda)^{\oplus b_\lambda}\subset (S^\lambda)^{\oplus a_\lambda}$ such that $I_d=\bigoplus_{\lambda\vdash dc} (S^\lambda)^{\oplus b_\lambda}$. This is equivalent to finding a $b_\lambda$-dimensional linear subspace in the space of highest weight vectors in $(S^\lambda)^{\oplus a_\lambda}$. We provide a database of polynomials in
$$S^d(S^c) := S^d(S^c((\CC^n)^*))$$
that for each $\lambda$ provides a basis of the highest weight space of $(S^\lambda)^{\oplus a_\lambda}$. Finally, we apply exact linear algebra methods to find which combinations of those vectors vanish on $X$. This is done by finding exact random points of $X$ giving linear conditions on highest weight spaces.

To generate a basis of the highest weight vectors in $S^d(S^c)$ one may first generate a basis of highest weight vectors of weight $\lambda$ in $(S^c)^{\otimes d}$. This is obtained by applying the Pieri rule. As writing this basis in terms of tensors is quite memory and time consuming, it is much better to simply remember it in terms of semistandard Young tableaux.
The symmetrizing operator $(S^c)^{\otimes d}\rightarrow S^d(S^c)$ maps this basis to a generating set. Out of that set one chooses a basis, using linear algebra. There are many choices to pick a basis out of a generating set. We choose a random initial element in the generating set and add it to our basis. Then we choose another random element in the generating set, and check if it is linearly independent to the current basis. If it is we add this new element to the basis. We repeat this process until the number of basis elements equals the multiplicity of $S^\lambda$ in $S^d(S^c)$. To check linear independence it is enough to be able to evaluate a polynomial corresponding to a given Young tableaux at many points. We apply a method that allows fast evaluation, without the necessity to expand the whole highest weight vector.

\section{Algorithm}\label{sec:algorithm}

We describe here how to convert a Young tableau into a highest weight polynomial over the monomial basis. Evaluation at a point in $X$ is then straightforward.
In this way, if we can sample efficiently from $X$, we can evaluate the basis of $a_\lambda$ many highest weight polynomials at $a_\lambda$ sampled points, obtain a square matrix $A$ of evaluations, and use linear algebra to compute $b_\lambda = \dim \text{ker} A$.

We are given two natural numbers $d,c \in \NN$.
Moreover, we are given a so-called \emph{isobaric Young tableau}. This is a top-left justified arrangement of $dc$ many boxes with entries from $\{1,\ldots,d\}$ such that every entry appears exactly $c$ times, see this example with $d=4, c=3$:
\ytableausetup{boxsize=1.1em}
\[
\ytableaushort{
111234,2234,34
% {1^2_1}{1^2_1}{1^2_1}{1^2_1}{1^2_1}{1^2_1}
% {1^2_2}{1^2_2}{1^2_2}{1^2_2}{1^2_2}{1^2_2}
% ,
% {*(lightgray)2_1}{*(lightgray)2_1}{*(lightgray)2_1}{*(lightgray)2_1}{*(lightgray)2_2}{*(lightgray)2_3}{*(lightgray)2_3}{*(lightgray)2_4}{*(lightgray)2_5}{*(lightgray)2_6}{*(lightgray)2_7}{*(lightgray)2_8}
% ,
% {3^2_1}{3^2_1}{3^2_1}{3^2_1}{3^2_1}{3^2_1}
% {3^2_2}{3^2_2}{3^2_2}{3^2_2}{3^2_2}{3^2_2}
% ,
% {4^2_1}{4^2_1}{4^2_1}{4^2_1}{4^2_1}{4^2_1}
% {4^2_2}{4^2_2}{4^2_2}{4^2_2}{4^2_2}{4^2_2}
}
%\young(111234,2234,34)
\]
In fact, we may assume that the tableau is \emph{semistandard}, which means that the entries are increasing within each column from top to bottom and they are nondecreasing within each row from left to right. The example above is semistandard.
It is an open question whether or not it is possible to use semistandardness to get a speed-up in the running time, see \cite{BDI:21} for the details.

We color the boxes in the same color if and only if they have the same number, and then we remove the numbers:
\[
\ytableaushort{
{*(red) }{*(red) }{*(red) }{*(blue) }{*(green) }{*(yellow) },{*(blue) }{*(blue) }{*(green) }{*(yellow) },{*(green) }{*(yellow) }
}
\]
Let $\mu_i$ denote the number of boxes in column $i$. Let $\aS_k$ denote the symmetric group on $k$ letters. A \emph{column permutation assignment} is defined as an assignment of numbers to the boxes such that in each column $i$ each number from $\{1,\ldots,\mu_i\}$ appears exactly once. For example, this is a column permutation assignment:
\[
\ytableaushort{
{*(red)1}{*(red)2}{*(red)2}{*(blue)1}{*(green)1}{*(yellow)1},{*(blue)2}{*(blue)1}{*(green)1}{*(yellow)2},{*(green)3}{*(yellow)3}
}
\]
Each column in a column permutation assignment specifies a permutation, so we can define the \emph{sign} of a column permutation assignment to be the product of the signs of the permutations that correspond to the columns. The example above has sign $1 \cdot (-1) \cdot (-1) \cdot 1 \cdot 1 \cdot 1 = 1$.

To each column permutation assignment $T$ we assign the word $w(T)$ that is obtained by reading from $T$ first all entries from one color, then from the next, and so on. The order of colors and the order in which we read entries from the same color does not matter, because we define two words of length $cd$ to be \emph{equivalent} if they arise from each other by permuting symbols within the block $\{1,\ldots,c\}$ or within $\{c+1,\ldots,2c\}$, and so on, or if they arise by permuting the $d$ many blocks (in other words, they are equivalent if and only if they lie in the same orbit under the action of the wreath product $\aS_c \wr \aS_d$). The equivalence class of words $w(T)$ in the example above is $\{\{1,2,2\},\{1,1,2\},\{1,1,3\},\{1,2,3\}\}$.
To every column permutation assignment~$T$, let $\kappa(T)$ denote the equivalence class of $w(T)$.
Consider the vector space spanned by all possible $\kappa(T)$, where we interpret distinct $\kappa(T)$ to be linearly independent unit vectors.

Up to a simple rescaling of the basis, the highest weight polynomial in monomial presentation is then
\begin{equation}\tag{$\dagger$}\label{eq:bruteforce}
\sum_{\text{column permutation assignment } T} \text{sgn}(T) \kappa(T).
\end{equation}
For example, let $c=d=2$ and take the tableau $\vcenter{\hbox{$ \ytableausetup{smalltableaux}
\ytableaushort{
{*(red) }{*(red) },{*(blue) }{*(blue) }
}
$}}$, then the set of column permutation assignments (written into the tableau) is
$\Big\{\vcenter{\hbox{$ \ytableausetup{smalltableaux}
\ytableaushort{
{*(red) 1}{*(red) 1},{*(blue) 2}{*(blue) 2}
}
$}}
,
\vcenter{\hbox{$ \ytableausetup{smalltableaux}
\ytableaushort{
{*(red) 1}{*(red) 2},{*(blue) 2}{*(blue) 1}
}
$}}
,
\vcenter{\hbox{$ \ytableausetup{smalltableaux}
\ytableaushort{
{*(red) 2}{*(red) 1},{*(blue) 1}{*(blue) 2}
}
$}}
,
\vcenter{\hbox{$ \ytableausetup{smalltableaux}
\ytableaushort{
{*(red) 2}{*(red) 2},{*(blue) 1}{*(blue) 1}
}
$}}
\Big\}
$,
so the sum \eqref{eq:bruteforce} becomes $2\{\{1,1\},\{2,2\}\}-2\{\{1,2\},\{1,2\}\}$. This corresponds to the tensor $2 (e_1\odot e_1) \odot (e_2\odot e_2) - 2 (e_1 \odot e_2)\odot (e_1 \odot e_2)$ in $S^2(S^2)$, where $a \odot b := \frac 1 2 (a \otimes b + b \otimes a)$.

We compute the sum \eqref{eq:bruteforce} in a brute force way and store the result in a file. This means that the file then contains the highest weight polynomial in the monomial basis. In particular, the evaluation of a highest weight polynomial from a file at any point is very efficient.
A bottleneck in the computation \eqref{eq:bruteforce} is the number of column permutation assignments. For example, if $d=8$, $c=4$, then for the Young diagram with row lengths $\lambda=(8,8,8,8)$ we have 110\,075\,314\,176 many column permutation assignments.
Therefore it is imperative to perform as few operations as possible for each summand. Here are a few points which accelerate the computation:
\begin{enumerate}
\item We use a Gray code to iterate through the sum so that the sign alternates for every summand. A Gray code is a way to iterate over the set of all permutations so that each permutation differs from the next by only a transposition.
Therefore in each step the sign of the permutation flips and hence we never have to compute the sign of a permutation.
The Gray code we use for each column is \emph{Algorithm P} in \cite[Sec.~7.2.1.2.]{Knuth} and we hardcode its list of permutations for each column.
\item We do not compute $\kappa(T)$, because it would require sorting a list of lists. Instead we use an \emph{equivariant hash function},
which is a function that takes list of $d$ lists of numbers that are each of length $c$ and assigns this list a number (its so-called \emph{hash value}) such that every reordering under the wreath product $\aS_c\wr\aS_d$ has the same hash value and in a way that lists of lists that are \emph{not} equivalent under the wreath product action get different hash values (this last property is called \emph{collision-freeness}).
We can efficiently compute the hash value for a column permutation assignment and the equivariance of the hash function guarantees that words that are equivalent under the wreath product action are mapped to the same hash value.
The collision-free hash function is chosen in a precomputation step.
\item To crucially speed up to computation
the hash value is not computed for each summand, but the hash value is just \emph{adjusted} at each step. This is possible, because the hash function is chosen as follows. Let $T_{i,j}$ be the $j$th entry in the $i$th colored block of $T$. Then, the hash function $h$ is
\[
h(T) := \sum_{i=1}^d \big(\sum_{j=1}^c \iota(T_{i,j})\big)^k \text{ mod }p
\]
for a suitable $k \in \NN$ and prime $p$,
where $\iota(i)$ is the $i$-th entry in a fixed array of random numbers from $\{0,\ldots,p-1\}$.
Raising to the $k$-th power is done by repeated squaring.
The Gray code ensures that only two blocks are changed and only one entry in each block, which makes updating the hash value very efficient. To give a rough idea of the performance, after the precomputation of the hash function the summation over the 110\,075\,314\,176 entries for $\lambda=(8,8,8,8)$ takes only a few hours on a laptop.
\end{enumerate}
Those ideas are incorporated into our implementation:
\begin{samepage}
\begin{algorithmic}
\Function{WriteHighestWeightPolyToFile}{isobaric Young tableau $\mathscr T$}
\State Repeatedly choose $k$ and $p$ until the hash function $h$ is collision free
\State Initialize an array $A$ of $p$ many integers
\State Initialize the hash value $\gamma$ for the first summand in the sum \eqref{eq:bruteforce}
\State $\alpha := 1$
\For {column permutation assignment $T$}
 \State $A[\gamma] := A[\gamma] + \alpha$
 \State efficiently update $\gamma$ to be the hash value of the next summand
 \State $\alpha := -\alpha$
\EndFor
\State Start a new empty highest weight polynomial file
\For {basis vector $v$ in $S^d(S^c)$}
\State Set $\beta := A[h(v)]$
\State Append ``+'' and $\beta$ and ``$\cdot$'' and $v$ to the highest weight polynomial file
\EndFor
\EndFunction
\end{algorithmic}
\end{samepage}
This works well as long as $p$ many integers can be stored in the memory and no hash collisions appear.
For extremely large cases this is a problem, but then we just accept some hash collisions are store them, and whenever a hash collision appears the values are hashed again with a second hash function. This is also done with more than 2 hash functions for extremely large problems.
The number of such hash functions is determined before running the algorithm and an estimate of the number of hash functions is obtained based on estimating the number of hash collisions using the birthday problem formula.

\section{Numerical methods}\label{sec:exm}
The algorithm that we have described in the last section is symbolic. It is based on exact computations, thus yielding exact results. As we have demonstrated by obtaining Theorem \ref{thm:main}, it can go beyond the cases that the classical method relying on Gr\"obner basis \cite{cox2013ideals, MiSt} can cope with.

Nevertheless, there are still limits to our algorithm with the current technology that numerical methods can surpass. For instance, our main theorem (Theorem~\ref{thm:main}) shows that no equations of degree at most 8 vanish on the variety of quartic symmetroids $Q_4$. But we could not find the minimal degree $d$, for which there are equations; i.e., such that $I_d\neq \emptyset$. Numerical methods, although not exact, can help to make an educated guess for those numbers. In this last section we want to explain this.

We first explain an approach on how to compute the degree of $X$.
Thereafter, we will discuss that one can in principle extract the minimal $d$, such that $I_d\neq \emptyset$, from this computation. This poses new numerical challenges, however.

\subsection{Cubic symmetroids}
Our algorithm was successful in case of the variety of cubic symmetroids $Q_3$. First we verified that there are no equations up to degree 10. This can be performed without problems on a laptop. In degree~11, the longest part is to transform the database of highest weight vectors from tableau to polynomials. This takes a few days. However, this should be considered as a precomputation and this database, once created, can be used in future for any other problem concerning $S^{11}(S^3(\CC^4))$. In particular, in this step we create six linearly independent polynomials of weight $[15,6,6,6]$.
Once this is done, our algorithm finds fast (within hours) the unique highest weight polynomial of degree $11$ in the ideal. This is a unique linear combination of the six highest weight vectors of weight $[15,6,6,6]$ that vanishes on $Q_3$.
We do not present this polynomial in the article, as it is quite large. It may be downloaded from \cite{url}.
Using the same algorithm we also determined all of the degree $12$ polynomials that appear in the ideal of $Q_3$. The obtained representations corresponded exactly to those that appear in the tensor product $[15,6,6,6]\otimes [3]$ (and correspond to partitions with at most four entries). This motivates the following definition and question.
\begin{definition}
A variety $X\subset \CC^n$ with a $G$ action is called \emph{$G$-principle} if the ideal $I(X)$ is generated by one, irreducible representation of $G$.
\end{definition}
We note that to provide the whole ideal of a $G$-principle variety, it is enough to know the group action and provide only one polynomial that generates the irreducible representation.
\begin{question}
Is the variety $Q_3$ $GL(4)$-principle?
\end{question}
The highest weight polynomials in degree $12$ may also be downloaded from \cite{url}.

\subsection{A numerical approach for computing the degree}
We make a numerical computation to determine the degree of the symmetroid $Q_3$ ~(\ref{def_Q}). The approach described in this section can easily be generalized to the general situation involving~$X$, but for simplicity we will stick to the special situation with~$Q_3$.

We use coordinates by setting the first matrix to be $A_0 = \mathrm{diag}(a_1,a_2,a_3)$, where $a_1,a_2,a_3$ are variables. Then, we have the following situation:
\begin{align*}
f:V&\to W,\\
 a=(a_1,a_2,a_3,A_1,A_2,A_3) &\mapsto \text{coefficients of }\det(x_0 A_0 + x_1A_1 + x_2A_2 + x_3A_3)
\end{align*}
and $\dim V = 21$ and $\dim W = 20$. Here, coefficients means the coefficients of a polynomial in $x$. To determine the dimension of $Q_3$ we compute the rank of the Jacobian matrix of $f$ at a random point. We get that this rank is $16$  (see~\cite{symmetroids2021}). Therefore, $\dim Q_3=16$.
This implies that the dimension of the fibers of~$f$ for a general point $h\in Q_3$ is
$\dim f^{-1}(h) = 5$.
The degree of $Q_3$ multiplied by the degree of a general fiber $f^{-1}(h), h\in Q_3,$ is thus the number of isolated complex solutions of the following system of $16$ polynomial equations in the $16$ variables~$b=(b_1,\ldots,b_{16})$:
\begin{equation}\label{system1}R \cdot f(a) = r \quad\text{and}\quad a = S\cdot b+s,\end{equation}
where
$R \in \mathbb{C}^{16 \times 20}$, $r \in \mathbb{C}^{16}$,
$S \in \mathbb{C}^{21\times 16}$ and $s \in \mathbb{C}^{21}$
are chosen randomly. We can exploit monodromy by varying the coefficients of $R$ and $r$ in loops and numerically tracking the solutions along those loops. This produces new solutions for~(\ref{system1}). The details of this technique are, for instance, explained in~\cite{DHJLLS2018}. An initial solution  for this system can easily be generated. In the computation it is enough to keep one point in each fiber $f^{-1}(h), h\in Q_3$, so that we do not have to compute the degree of a general fiber to get the degree of~$Q_3$. A concrete implementation is at \cite{symmetroids2021}.

We can use the same method for the quartic symmetroid $Q_4$. In this case we have
$\dim V = 34$, $\dim W = 35$, and $\dim Q=25$. The dimension of the fibers of $f$ for a general point is
$9$. However, we could not finish this computation. We stopped the computation manually. At this point $849998$ solutions for (\ref{system1}) had been found.

\subsection{The degree of $Q_3$}
We used \texttt{HomotopyContinuation.jl} \cite{BT2018} for the computation of degrees of symmetroids. The code for the cubic symmetroid can be found at \cite{symmetroids2021}. We compute $305$ solutions of the system~(\ref{system1}) for $Q_3$. On a laptop this takes about two minutes. We use the certification method \cite{breiding2020certifying} implemented in \texttt{HomotopyContinuation.jl} to show that the 305 solutions we have found correspond to 305 true distinct solutions for~(\ref{system1}). This confirms that the degree of $Q_3$ is at least 305. It was proven in~\cite{Vainsencher2003} that the degree is equal to 305.

%For the quartic symmetroid we executed a three months  computation, after which the algorithm had found $849998$ solutions for (\ref{system1}). At this point the computation was aborted manually, because it hadn't found any more solution in a week. This led us to state Conjecture \ref{conj1}.

\subsection{Further directions}
Here, we explain an approach for answering the following question: Given a homogeneous polynomial map $f:\CC^a\rightarrow \CC^b$ what is the dimension of the vector space $I_d$ of polynomials of degree $d$ that vanish on the image?

The basic idea is this: suppose that we have run the algorithm from the previous section. Then, we have found a linear space $L$ in $W=S^c(\CC^n)$ and points $w_1,\ldots,w_\delta \in X\cap L$, such that $\delta$ is the degree of $X$. Any equation that vanishes on~$X$ also vanishes on the $w_i$. We now discuss when the reverse is true. If this holds, we can check numerically by solving a system of linear equations, whether or not there are equations of a fixed degree $d$ vanishing on the $X$. Note that this does \emph{not} yield equations for $X$. Furthermore, we can use coordinates for $L$ for doing the linear algebra. This kind of \emph{dimensionality reduction} can provide a significant reduction in computational complexity. These ideas first appeared in \cite{HIL:13}.

Let us write $b:=\dim(S^c(\CC^n))$ and the image of $f$ is invariant under the action of $\GL(\CC^n)$. We ask for the dimension of $I_d$, that is the degree $d$ part of $I$. It should be emphasized that this is naturally a problem in linear algebra, as $I_d$ is a vector space. Each point $x\in X$ determines a linear condition on the space $S^d(\CC^b)$, giving rise to a hyperplane containing $I_d$. In fact, $I_d$ is the intersection of all such hyperplanes. For dimensional reasons, it would be enough to pick consecutively random $x\in X$ and intersect the hyperplanes in $S^d(\CC^b)$, until the intersection stabilizes. This is indeed sometimes done in practice, but the main problem is the large dimension ${d+b-1}\choose{d}$ of the space $S^d(\CC^b)$. The method we describe is particularly useful if:
\begin{enumerate}
\item the codimension of $X$ is small,
\item the degree of $X$ is small.
\end{enumerate}
From now on we work in the projective space $\PP(\CC^b)$ and consider $X$ as a projective variety. Let $e:=\codim X$.

We pick a random subspace $L=\PP^e\subset\PP(\CC^b)$. By Bertini's theorem $\PP^e$ intersects $X$ in $\delta=\deg X$ many smooth points
$$S=L\cap X.$$
A random linear form $h_1$ is not a zero divisor in the ring $\CC[y_1,\dots,y_b]/I$, hence we have an exact sequence:
$$0\rightarrow \CC[y_1,\dots,y_b]/I\rightarrow \CC[y_1,\dots,y_b]/I\rightarrow \CC[y_1,\dots,y_b]/(I+(h_1))\rightarrow 0,$$
where the first map is multiplication by $h_1$. Hence, the Hilbert series of $I+(h_1)$ equals $(1-t)$ times the Hilbert series of $I$. In particular, the numerators of the Hilbert series are the same. The number of linear $h_1,\dots,h_l$ such that $h_{i+1}$ is not a zero divisor modulo $I+(h_1,\dots,h_i)$ for every $0\leq i<l$ is governed by the depth of the (localization of the) ring $\CC[y_1,\dots,y_b]/I$. Depth is always at most equal to the dimension and the cases when equality holds are called Cohen-Macaulay.

After choosing $e=\codim X$ many linear forms, we arrive at the ring
$$\CC[y_1,\dots,y_b]/(I+(h_1,\dots h_e))$$
which describes $S$ as a \emph{projective} scheme. In general, the ideal $(I+(h_1,\dots h_e))$ may have an embedded component at zero, however this again does not happen if $X$ is arithmetically Cohen-Macaulay (which means that its coordinate ring is Cohen-Macaulay). In practice, we next choose an affine linear form and add it to the ideal to represent $S$ as a finite subset of an affine space.

In particular, if our variety $X$ is arithmetically Cohen-Macaulay
then the Hilbert function of the finite set $S$ in fact encodes the numerator of the Hilbert series of~$X$. In any case a nonzero element in $I_d$ gives rise to a nonzero element in $I(S)_d$. As long as $I(S)_d=0$ we also have $I_d=0$, hence we do not have to look for equations in those degrees. Further, in smallest degree $d$ such that $I(S)_d\neq 0$ we have $\dim I_d= \dim I(S)_d$ in the Cohen-Macaulay case.

\begin{ex}
In the following example we construct a toric ring of small depth. Consider the map for $x=(x_1,x_2,x_3,x_4)$:
$$x\mapsto (x_1x_2^4,x_1x_2^3x_3,x_1x_2x_3^3,x_1x_3^4,x_1x_2^4x_4,x_1x_2^3x_3x_4,x_1x_2x_3^3x_4,x_1x_3^4x_4).$$
The image is a toric variety of projective dimension two and degree eight. It is minimally generated by nine quadrics and twelve cubics. If we intersect the image with two affine linear forms we obtain eight points. These eight points do not contribute to new linear equations, however their ideal has thirteen minimal generators in degree two.
\end{ex}

Thus, if we know $S$, we may estimate $\dim I_d$ using linear algebra approach described above, but now we deal with points in the ambient space of dimension $e=\codim X$. Hence, we have to solve $\deg X$ many \emph{linear} equations in ${d+e-1}\choose{d}$ many variables.

Numerical methods help us both: to obtain $S$ and to solve the linear equations. To generate $\PP^e$ we take a span of $e+1$ many random/general points of $X$. We obtain~$\PP^e$ together with ${e+1}$ many points of $S$. To generate  all of $S=\{w_1,\ldots,w_\delta\}$ we apply the monodromy method from the previous subsection.

Now a new problem arises. As our points are just approximations of the points in $S$, if we ask for the rank of the matrix associated to the system of linear equations, symbolically it will always be nondegenerate. Further, the matrix we obtain depends on the choice of the basis of degree $d$ polynomials we take. The idea is to look at the singular values of the associated matrix in the basis. This allows us to discover the rank of the approximated matrix.

We plan to apply the approach, that we have just described, to computations involving $\mathrm{GL}$ invariant families of polynomials. %The particular case of quartic symmetroids $Q_4$ is the next challenge.
%, we were not able to apply it, yet. We first must settle Conjecture \ref{conj1}.

\bibliography{bibML}
\bibliographystyle{plain}
\end{document}